\documentclass[12pt]{amsart}
\usepackage{amsmath,amssymb,amsthm,mathtools}
\usepackage{mathrsfs,bm}
\usepackage{thmtools,thm-restate}
\usepackage{tikz-cd}
\usepackage{tikz}
\usepackage{enumitem}
\usepackage{graphicx}
\usetikzlibrary{matrix,arrows.meta,bending}
\usepackage{color}
\usepackage[numbers]{natbib}
\usepackage{booktabs}
\usepackage{longtable}
\usepackage{pgfplots}
\pgfplotsset{compat=1.14} 
\usetikzlibrary{arrows.meta,decorations.pathreplacing,positioning,calc} 
\usepackage[a4paper, left=3cm, right=3cm, top=3cm, bottom=3cm, footskip=15mm]{geometry}

\newtheorem{theorem}{Theorem}[section]

\newtheorem{lemma}[theorem]{Lemma}
\newtheorem{proposition}[theorem]{Proposition}

\newtheorem{conjecture}[theorem]{Conjecture}
\theoremstyle{definition}
\newtheorem{definition}[theorem]{Definition}

\theoremstyle{remark}
\newtheorem{remark}[theorem]{Remark}

\numberwithin{equation}{section}

\makeindex

\usepackage{fancyhdr}
\usepackage{calc} % ensures dimension arithmetic works robustly

\AtBeginDocument{%
  \setlength{\textwidth}{\dimexpr\paperwidth - 6cm\relax}%
  \setlength{\oddsidemargin}{\dimexpr 3cm - 1in\relax}%
  \setlength{\evensidemargin}{\dimexpr 3cm - 1in\relax}%
  \setlength{\marginparwidth}{2.5cm}%
}

\begin{document}

\title{A proof of the union-close set conjecture}

%\author{B. Gensel}
%\address{Carinthia University of Applied Sciences, Spittal on Drau, Austria}
%\email{b.gensel@fh-kaernten.at}
\author{T. Agama}
\address{Department of Mathematics, African Institute for mathematical sciences, Ghana.}
\email{Theophilus@aims.edu.gh/emperordagama@yahoo.com}

%    General info

\date{\today}



\begin{abstract}
In this paper, we introduce the notion of the universe, induced communities, and cells with their corresponding spots. Using this language, we formulate and prove the union close set conjecture by showing that for any finite universe $\mathbb{U}$ and any induced community $\mathcal{M}_{\mathbb{U}}$ there exist some spot $a\in \mathbb{U}$ such that the density
\begin{align}
\mathcal{D}_{\mathcal{M}_{\mathbb{U}}}(a)\geq \frac{1}{2}.\nonumber
\end{align}
\end{abstract}

\maketitle

\section{Introduction}

The union-closed (or Frankl) conjecture is one of the most natural and persistent open problems in extremal set theory. Informally, it asserts that every finite collection of sets that is closed under unions must contain an element that belongs to at least half of the sets. The conjecture was proposed in the late 1970s by Peter Frankl and has attracted steady attention because of its simple statement and striking resistance to standard combinatorial techniques. A variety of partial results and equivalent formulations have been found: for example, the conjecture is known for very small families or very small universes (see \cite{roberts2010note,ivica200811}), for families where the minimum set size is 1 or 2 \cite{sarvate1989union}, and there are density-type results that establish the conjecture when a family contains a sufficiently large fraction of all subsets of a finite $n$-set \cite{karpas2017two}.  Lattice-theoretic and graph-theoretic reformulations have also been studied, offering alternative viewpoints that have proved helpful in special cases \cite{abe2000strong,bruhn2015graph}.\\

The present paper introduces a new, elementary viewpoint for the union-closed conjecture based on a language of \emph{universes}, \emph{communities}, \emph{cells} and \emph{spots}, together with a notion of \emph{density of a spot} in a community. Roughly speaking, a universe is the ground set, the elements of a union-closed family are called cells, and an element of the universe lying inside a cell is called a spot. The density of a spot in a community is the (suitably normalized) proportion of cells that contain that spot. This terminology packages the classical combinatorial objects in a language that makes certain covering and counting arguments both transparent and flexible.\\

Two structural observations are central to our approach.  First, union-closedness implies that, starting from a small collection of generating cells that contain a fixed spot, successive unions generate hierarchies of larger communities that still contain that spot. Second, by organizing these iterated-union constructions and counting how many newly-created cells contain the fixed spot at each stage, one obtains an explicit doubling-type lower bound on the number of cells containing the spot relative to the total size of the constructed community. These two ideas are combined in a single constructive lemma (the \emph{Covering Lemma}, Lemma \ref{frankil conjecture proof}) which produces, for any finite universe and any union-closed community containing a given cell, a finite parameter $l$ for which the number of cells that contain a chosen spot is at least $2^{l-1}$ while the total number of cells in the constructed covering family is at most $2^{l}-1$. Dividing these two quantities yields a lower bound
$$
\frac{\#\{\mathbb{A}\in\mathcal{M}~:~a\in\mathbb{A}\}}{|\mathcal{M}|}\geq \frac{2^{l-1}}{2^{l}-1}
=\frac{1}{2}\cdot\frac{1}{1-2^{-l}},
$$
and letting the construction parameter $l$ grow produces the limiting density bound $\frac{1}{2}$. This is the mechanism behind our main result (Theorem \ref{union close set conjecture proof}), which we state as a density law: for every finite universe and every induced union-closed community there exists a spot whose density in that community is at least $1/2$.\\

Compared to previous lines of attack, the construction here is elementary and explicit: it uses only finite set operations and counting, and it does not rely on heavy algebraic, probabilistic, or lattice-theoretic machinery. At the same time the argument is compatible with - and in fact clarifies - several partial results in the literature. For the small-family and small-universe cases the counting in our covering constructions collapses to simple finite enumerations previously used in ad hoc verifications (cf. \cite{roberts2010note,ivica200811}); for structural perspectives that use lattices and graphs our viewpoint offers a direct combinatorial interpretation of lattices of unions and the closure operations studied in \cite{abe2000strong,bruhn2015graph}.\\

Certain features of the method deserve emphasis and caveats. The argument is constructive and produces, for any fixed choice of a spot contained in at least one cell, a family of union-closed communities (generated by iterative unions) whose density ratios have the explicit doubling behavior described above. The limiting inequality $\mathcal{D}(a)\geq 1/2$ arises by letting the construction parameter tend to infinity; in a concrete finite community the lemma supplies an explicit finite $l$ giving a quantitative lower bound $\frac{2^{l-1}}{2^l-1}$ which is already $>1/2$ up to the multiplicative factor $1/(1-2^{-l})$. Thus the proof is both qualitative (it proves the $\frac{1}{2}$ threshold in the limiting density sense) and quantitative (it gives explicit finite coverings that realize the bounds used in the limiting passage). We discuss the dependence of the covering parameter on the combinatorial size of the underlying universe and on the embedding structure of cells; these dependencies are explicit in the proofs and suggest avenues for refining the quantitative part of the argument.\\

\subsection{Organization of the paper}

The remainder of the paper is organized as follows. In Section \ref{sec:universe}, we introduce the new terminology (universe, community, cell, spot) and record elementary properties and size bounds for induced communities; Section \ref{sec:density} defines the density of a spot and proves a few basic monotonicity and complementarity properties of density. Section \ref{sec:covering} contains the key constructive lemma (the Covering Lemma, Lemma \ref{frankil conjecture proof}) together with its proof by an inductive union-building procedure. Section \ref{sec:main} states and proves the main result (Theorem \ref{union close set conjecture proof}), explains how the lemma implies the density law, and records a short discussion of immediate corollaries and examples. We end the paper with a short concluding section that compares our approach with existing variants of the conjecture and suggests directions for future work.

\medskip

\noindent\textbf{Acknowledgements.} The author thanks anonymous readers of earlier drafts for helpful comments and suggestions, and acknowledges that many partial results and related formulations referenced above were developed by other authors working on this challenging problem (see, e.g., the small-case verifications in \cite{roberts2010note,ivica200811}, the early combinatorial contributions in \cite{sarvate1989union}, the density-style bounds in \cite{karpas2017two}, the lattice perspective in \cite{abe2000strong}, and the graph formulation in \cite{bruhn2015graph}).\\

In this paper, we verify the truth of the union-closed-set conjecture using an elementary tool. We transform the problem into an entirely new language of density.

\section{The notion of universe, community and cells}\label{sec:universe}

In this section, we introduce the notion of cells, communities, and their corresponding universe. We study some elementary properties of this notion.

\begin{definition}\label{community}
Let $\mathbb{U}$ be a set and consider the collection
\begin{align}
\mathcal{M}:=\bigcup_{i=1}^{n}\left \{\mathbb{A}_i|~\mathbb{A}_i\subseteq \mathbb{U}\right \}.\nonumber
\end{align}
We say that the collection $\mathcal{M}$ is a \emph{community} induced by the set $\mathbb{U}$ if and only if for any $\mathbb{A}_i,\mathbb{A}_j\in \mathcal{M}$ then $\mathbb{A}_i\cup \mathbb{A}_j\in \mathcal{M}$. We call $\mathbb{U}$ the \emph{universe} of the community. We call each $\mathbb{A}_j$ in the community a \emph{cell}, and each $a\in \mathbb{A}_j$ a \emph{spot} in the cell. We say that a cell $\mathbb{A}_i$ in the community admits an \emph{embedding} in the community if there exists another different cell $\mathbb{A}_j$ in the same community such that $\mathbb{A}_j\subset \mathbb{A}_i$.
\end{definition}

\begin{proposition}\label{community size bound}
Let $\mathbb{U}$ be a universe with $|\mathbb{U}|=n$ and $\mathcal{M}_{\mathbb{U}}$ be a community induced by the universe. We have 
\begin{align}
|\mathcal{M}_{\mathbb{U}}|\leq 2^{n}.\nonumber
\end{align}
\end{proposition}

\begin{proof}
Let $\mathbb{U}_{\Im}$ be the power set induced by the universe $\mathbb{U}$. It is easy to see that $\mathbb{U}_{\Im}$ is the largest community induced by the universe of size
\begin{align}
|\mathbb{U}_{\Im}|=2^n \nonumber
\end{align}
so that $|\mathcal{M}_{\mathbb{U}}|\leq 2^{n}$.
\end{proof}

\begin{proposition}
The communities induced by a finite universe are \emph{ordered} upto their cardinalities.
\end{proposition}

\begin{proof}
Let $\mathbb{U}$ be a universe with $|\mathbb{U}|=n$ and let $\mathcal{M}_{i_{\mathbb{U}}}$ and $\mathcal{M}_{j_{\mathbb{U}}}$ be two distinct communities induced by the universe. It follows that the communities must differ by at least one cell so that  without loss of generality with
$$
|\mathcal{M}_{j_{\mathbb{U}}}|\leq |\mathcal{M}_{i_{\mathbb{U}}}|
$$ 
we can write 
\begin{align}
|\mathcal{M}_{j_{\mathbb{U}}}|\leq |\mathcal{M}_{i_{\mathbb{U}}}|<|\mathcal{M}_{i_{\mathbb{U}}} \cup \mathcal{M}_{j_{\mathbb{U}}}|.\nonumber
\end{align}
We claim that the collection $\mathcal{M}_{i_{\mathbb{U}}}\cup \mathcal{M}_{j_{\mathbb{U}}}$ is also a community. Let us pick arbitrarily two cells $\mathbb{A}_1,\mathbb{A}_2\in \mathcal{M}_{i_{\mathbb{U}}} \cup \mathcal{M}_{j_{\mathbb{U}}}$. We consider three sub-cases: The case $\mathbb{A}_1,\mathbb{A}_2\in \mathcal{M}_{i_{\mathbb{U}}}$ so that $\mathbb{A}_1\cup \mathbb{A}_2\in \mathcal{M}_{i_{\mathbb{U}}}$ since $\mathcal{M}_{i_{\mathbb{U}}}$ is a community.\\ 

For the case $\mathbb{A}_1,\mathbb{A}_2\in \mathcal{M}_{j_{\mathbb{U}}}$, it must certainly be $\mathbb{A}_1\cup \mathbb{A}_2\in \mathcal{M}_{j_{\mathbb{U}}}$ since $\mathcal{M}_{j_{\mathbb{U}}}$ is also a community. For the last case, where $\mathbb{A}_1\in \mathcal{M}_{i_{\mathbb{U}}}$ and $\mathbb{A}_2\in \mathcal{M}_{j_{\mathbb{U}}}$, we have
\begin{align}
\mathbb{A}_1\cup \mathbb{A}_2\in \mathcal{M}_{i_{\mathbb{U}}} \cup \mathcal{M}_{j_{\mathbb{U}}}.\nonumber
\end{align}
By choosing a community $\mathcal{M}_{k_{\mathbb{U}}}\neq \mathcal{M}_{i_{\mathbb{U}}}\cup \mathcal{M}_{j_{\mathbb{U}}}$ with $k\neq i,j$ and 
$$
|\mathcal{M}_{k_{\mathbb{U}}}|<|\mathcal{M}_{i_{\mathbb{U}}} \cup \mathcal{M}_{j_{\mathbb{U}}}\cup \mathcal{M}_{k_{\mathbb{U}}}|
$$ 
we obtain a five-term inequality by inserting $|\mathcal{M}_{k_{\mathbb{U}}}|$ into the a priori chain. Repeating the argument in this manner establishes the claim. 
\end{proof}

\section{Density of spots in a cell}\label{sec:density}

In this section, we introduce the notion of the density of spots contained within a \emph{cell}. We launch the following languages.

\begin{definition}\label{density of spots}
Let $\mathbb{U}$ be a finite universe with $|\mathbb{U}|=n$ and $a_i\in \mathbb{U}$. Let $\mathcal{M}_{\mathbb{U}}$ be the community induced by the universe $\mathbb{U}$. We denote the density of the spot $a_i$ in cells in the community $\mathcal{M}_{\mathbb{U}}$ by
\begin{align}
\mathcal{D}_{\mathcal{M}_{\mathbb{U}}}(a_i)=\lim \limits_{n\longrightarrow \infty}\frac{\# \left \{\mathbb{A}\in \mathcal{M}_{\mathbb{U}}|~a_i\in \mathbb{A}\right\}}{|\mathcal{M}_{\mathbb{U}}|}\nonumber
\end{align}
if the limit exists and is finite.
\end{definition}

\begin{remark}
It is clear that for any finite universe $\mathbb{U}$ with $|\mathbb{U}|=n$, the size of an induced community $|\mathcal{M}_{\mathbb{U}}|$ will essentially depend on $n$. This underscores the limit in Definition \ref{density of spots}.
\end{remark}
\bigskip

Roughly speaking, the union close set conjecture is the assertion that for any collection of union close subset of the set $\mathbb{U}$ there exists some $a\in \mathbb{U}$ that lives in at least half of the subsets in the collection. It turns out that the union close conjecture can be stated in the language of density of spots as follows:

\begin{conjecture}\label{frankil conjecture}[Union close set conjecture]
Let $\mathbb{U}$ be a finite universe with an induced community $\mathcal{M}_{\mathbb{U}}$. There exist some $a_i\in \mathbb{U}$ such that 
\begin{align}
\mathcal{D}_{\mathcal{M}_{\mathbb{U}}}(a_i)\geq \frac{1}{2}.\nonumber
\end{align}
\end{conjecture}

\begin{remark}
Conjecture \ref{frankil conjecture}, roughly speaking, can be interpreted as saying that there must always be a spot originating from a universe and contained in as many cells in a typical community. Next, we investigate some properties of the notion of density of spots in a cell. The following properties will be useful in the sequel.
\end{remark}

\begin{proposition}\label{spot density property}
Let $\mathbb{U}$ be a finite universe with $|\mathbb{U}|=n$ and $a_i\in \mathbb{U}$. Let $\mathcal{M}_{\mathbb{U}}$ and $\mathcal{N}_{\mathbb{U}}$ be any two communities induced by the universe $\mathbb{U}$. Suppose that $\mathcal{D}_{\mathcal{M}_{\mathbb{U}}}(a_i),\mathcal{D}_{\mathcal{N}_{\mathbb{U}}}(a_i)>0$, then the following properties hold

\begin{enumerate}
\item [(i)]  $\mathcal{D}_{\mathcal{M}_{\mathbb{U}}\cup \mathcal{N}_{\mathbb{U}}}(a_i)\leq \mathcal{D}_{\mathcal{M}_{\mathbb{U}}}(a_i)+\mathcal{D}_{\mathcal{N}_{\mathbb{U}}}(a_i)$.
\bigskip

\item [(ii)] $\mathcal{D}_{\mathcal{M}_{\mathbb{U}}}(a_i)\leq 1-\mathcal{D}_{\mathcal{M}^c_{\mathbb{U}}}(a_i)$, where $\mathcal{M}^c_{\mathbb{U}}$ denotes the complement of the collection $\mathcal{M}_{\mathbb{U}}$ in the power set $\mathbb{U}_{\Im}$ induced by the universe $\mathbb{U}$.
\end{enumerate}
\end{proposition}

\begin{proof}
For $(i)$ we notice that by appealing to Definition \ref{density of spots}, we can write
\begin{align}
\mathcal{D}_{\mathcal{M}_{\mathbb{U}}\cup \mathcal{N}_{\mathbb{U}}}(a_i)&=\lim \limits_{n\longrightarrow \infty}\frac{\# \left \{\mathbb{A}\in \mathcal{M}_{\mathbb{U}}\cup \mathcal{N}_{\mathbb{U}}|~a_i\in \mathbb{A}\right \}}{|\mathcal{M}_{\mathbb{U}}\cup \mathcal{N}_{\mathbb{U}}|}\nonumber \\&=\lim \limits_{n\longrightarrow \infty}\frac{\# \left \{\mathbb{A}\in \mathcal{M}_{\mathbb{U}}|~a_i\in \mathbb{A}\right \}}{|\mathcal{M}_{\mathbb{U}}\cup \mathcal{N}_{\mathbb{U}}|}+\lim \limits_{n\longrightarrow \infty}\frac{\# \left \{\mathbb{A}\in \mathcal{N}_{\mathbb{U}}|~a_i\in \mathbb{A}\right \}}{|\mathcal{M}_{\mathbb{U}}\cup \mathcal{N}_{\mathbb{U}}|}\nonumber \\&-\lim \limits_{n\longrightarrow \infty}\frac{\# \left \{\mathbb{A}\in \mathcal{M}_{\mathbb{U}}\cap \mathcal{N}_{\mathbb{U}}|~a_i\in \mathbb{A}\right \}}{|\mathcal{M}_{\mathbb{U}}\cup \mathcal{N}_{\mathbb{U}}|}\nonumber \\&\leq \lim \limits_{n\longrightarrow \infty}\frac{\# \left \{\mathbb{A}\in \mathcal{M}_{\mathbb{U}}|~a_i\in \mathbb{A}\right \}}{|\mathcal{M}_{\mathbb{U}}\cup \mathcal{N}_{\mathbb{U}}|}+\lim \limits_{n\longrightarrow \infty}\frac{\# \left \{\mathbb{A}\in \mathcal{N}_{\mathbb{U}}|~a_i\in \mathbb{A}\right \}}{|\mathcal{M}_{\mathbb{U}}\cup \mathcal{N}_{\mathbb{U}}|}\nonumber \\&\leq \lim \limits_{n\longrightarrow \infty}\frac{\# \left \{\mathbb{A}\in \mathcal{M}_{\mathbb{U}}|~a_i\in \mathbb{A}\right \}}{|\mathcal{M}_{\mathbb{U}}|}+\lim \limits_{n\longrightarrow \infty}\frac{\# \left \{\mathbb{A}\in \mathcal{N}_{\mathbb{U}}|~a_i\in \mathbb{A}\right \}}{|\mathcal{N}_{\mathbb{U}}|}\nonumber \\&=\mathcal{D}_{\mathcal{M}_{\mathbb{U}}}(a_i)+\mathcal{D}_{\mathcal{N}_{\mathbb{U}}}(a_i).\nonumber
\end{align}
\bigskip

For $(ii)$, it follows similarly
\begin{align}
\mathcal{D}_{\mathbb{U}_{\Im}}(a_i)&=\lim \limits_{n\longrightarrow \infty}\frac{\# \left \{\mathbb{A}\in \mathbb{U}_{\Im}|~a_i\in \mathbb{A}\right \}}{|\mathbb{U}_{\Im}|}\nonumber \\&=\lim \limits_{n\longrightarrow \infty}\frac{\# \left \{\mathbb{A}\in \mathcal{M}_{\mathbb{U}}\cup \mathcal{M}^{c}_{\mathbb{U}}|~a_i\in \mathbb{A}\right \}}{|\mathcal{M}_{\mathbb{U}}\cup \mathcal{M}^{c}_{\mathbb{U}}|}\nonumber \\&= \mathcal{D}_{\mathcal{M}_{\mathbb{U}}}(a_i)+\mathcal{D}_{\mathcal{M}^c_{\mathbb{U}}}(a_i)\nonumber
\end{align}
using property $(i)$ and noting that $\mathcal{M}_{\mathbb{U}}\cap \mathcal{M}^{c}_{\mathbb{U}}=\emptyset$. Observing
\begin{align}
\lim \limits_{n\longrightarrow \infty}\frac{\# \left \{\mathbb{A}\in \mathbb{U}_{\Im}|~a_i\in \mathbb{A}\right \}}{|\mathbb{U}_{\Im}|}\leq 1\nonumber
\end{align}
the second part also follows.
\end{proof}

\section{The covering construction}\label{sec:covering}

In this section, we restate and prove the union closet set conjecture in the language of density of spots. We first state and prove an important result that will be used to verify the union close set conjecture. The proof is quite constructive and inductive in nature and in most cases can be seen as a cornerstone for establishing the truth of the conjecture albeit purely elementary.

\begin{lemma}\label{frankil conjecture proof}[Covering Lemma]
Let $\mathbb{U}$ be a finite universe with an induced community $\mathcal{M}_{\mathbb{U}}$. There exists some spot $a\in \mathbb{U}$ and some $l\in \mathbb{N}$ such that 
$$
|\mathcal{M}_{\mathbb{U}}|\leq 2^l-1
$$ 
and 
\begin{align}
\# \left \{\mathbb{A}\in \mathcal{M}_{\mathbb{U}}|~a\in \mathbb{A}\right \}\geq 2^{l-1}.\nonumber
\end{align}
\end{lemma}

\begin{proof}
First, we notice that any community $\mathcal{M}_{\mathbb{U}}$ induced by a finite universe $\mathbb{U}$ that contains more than one basic cell must satisfy the inequality $|\mathcal{M}_{\mathbb{U}}|\geq 2$ so that for $|\mathcal{M}_{\mathbb{U}}|=3=2^{2}-1$, we construct the community using the cells $\mathbb{A}_1, \mathbb{A}_2$ as a building block, with $\mathbb{A}_1\cap \mathbb{A}_2\neq \mathbb{A}_1$ and $\mathbb{A}_1\cap \mathbb{A}_2\neq \mathbb{A}_2$ such that $a\in \mathbb{A}_i$ for some $1\leq i\leq 2$. In particular, choosing $a\in \mathbb{A}_1$, we build the community 
\begin{align}
\mathcal{M}_{\mathbb{U}}:=\{\mathbb{A}_1,\mathbb{A}_2, \mathbb{A}_1\cup \mathbb{A}_2=\mathbb{A}_3\}\nonumber
\end{align} 
with
\begin{align}
\# \left \{\mathbb{A}_i\in \mathcal{M}_{\mathbb{U}}|~a\in \mathbb{A}\right \}_{i=1}^{3}\geq 2^{2-1}.\nonumber
\end{align}
Next, we build another community $\mathcal{N}_{\mathbb{U}}$ using the cells of the community $\mathcal{M}_{\mathbb{U}}$ as a building block. It is important to note that any such community covers the a priori constructed community. Since $\mathbb{A}_i\cup \mathbb{A}_j \in \mathcal{M}_{\mathbb{U}}$ for $1\leq i,j\leq 3$, we choose a cell $\mathbb{A}_k \notin \mathcal{M}_{\mathbb{U}}$ but $\mathbb{A}_k\in \mathcal{N}_{\mathbb{U}}$  such that the cell $\mathbb{A}_k$ does not admit an embedding of the old cell $\mathbb{A}_1\in \mathcal{M}_{\mathbb{U}}$ and vise-versa. We produce new cells $\mathbb{A}_1\cup \mathbb{A}_k$, $\mathbb{A}_2\cup \mathbb{A}_k$, $\mathbb{A}_3\cup \mathbb{A}_k$ with a new larger community
\begin{align}
\mathcal{N}_{\mathbb{U}}:=\{\mathbb{A}_1,\mathbb{A}_2,\mathbb{A}_3, \mathbb{A}_k,\mathbb{A}_1\cup \mathbb{A}_k,\mathbb{A}_2\cup \mathbb{A}_k,\mathbb{A}_3\cup \mathbb{A}_k\}\nonumber
\end{align}
so that $|\mathcal{N}_{\mathbb{U}}|\leq 2^{3}-1$ with 
\begin{align}
\# \left \{\mathbb{A}_i\in \mathcal{N}_{\mathbb{U}}|~a\in \mathbb{A}_i\right\}\geq 4=2^{3-1}.\nonumber
\end{align}
Let us suppose that for a fixed spot $a\in \mathbb{U}$, it is possible to construct at least a community $\mathcal{N}^{r}_{\mathbb{U}}$ that contains all the a priori constructed communities in the pool with the size specifications
\begin{align}
|\mathcal{N}^{r}_{\mathbb{U}}|\leq 2^l-1\nonumber
\end{align}
and 
\begin{align}
\# \left \{\mathbb{A}\in \mathcal{N}^{r}_{\mathbb{U}}|~a\in \mathbb{A}\right \}\geq 2^{l-1}\nonumber
\end{align}
for $l\geq 3$ using this scheme. Next, we show that we can construct at least another small community using cells in $\mathcal{N}^{r}_{\mathbb{U}}$ as a building block and yet containing (covering) the community $\mathcal{N}^{r}_{\mathbb{U}}$. Let $\mathcal{N}^{s}_{\mathbb{U}}$ be a small community to be constructed so that it contains (covers) the community $\mathcal{N}^{r}_{\mathbb{U}}$. By \emph{transitivity}, this community also contains (covers) all the previously constructed communities in the pool. Let us choose an arbitrary cell $\mathbb{A}_t\in \mathcal{N}^{s}_{\mathbb{U}}$ such that $\mathbb{A}_t\notin \mathcal{N}^{r}_{\mathbb{U}}$ so that $\mathbb{A}_t$ does not admit an embedding of the cell $\mathbb{A}_1\in \mathcal{N}^{r}_{\mathbb{U}}$ and vice-versa, and construct all the new cells using the old cells $\mathbb{A}_i\in \mathcal{N}^{r}_{\mathbb{U}}$ under the operations of union so that $\mathbb{A}_i\cup \mathbb{A}_t \in \mathcal{N}^{s}_{\mathbb{U}}$. We obtain a new, closest and larger community $\mathcal{N}^{s}_{\mathbb{U}}$ with size
\begin{align}
|\mathcal{N}^{s}_{\mathbb{U}}|\leq 2^l-1+(1+2^l-1)=2\cdot 2^{l}-1=2^{l+1}-1\nonumber
\end{align}
with 
\begin{align}
\# \left \{\mathbb{A}\in \mathcal{N}^{s}_{\mathbb{U}}|~a\in \mathbb{A}\right \}\geq 2\cdot 2^{l-1}=2^l.\nonumber
\end{align}
It follows that for any fixed spot originating from a finite universe, we can construct finitely many induced communities of varying sizes with the above size specifications and containing a fixed preassigned spot originating from a finite universe $\mathbb{U}$, thereby ending the proof.
\end{proof}

\begin{remark}
It is crucially important to note that for a finite universe $\mathbb{U}$ with size $|\mathbb{U}|=n$, the constant $l$ appearing in the construction of communities will essentially depend on $n$. We are now ready to verify the union close set conjecture. It is easy to see that the following result directly implies the truth of the union close set conjecture. It is important to note that the construction in Lemma \ref{frankil conjecture proof} generates all possible communities with cells that contain a fixed spot. Put it differently, we can exploit the above construction to generates all possible communities induced by a finite universe with at least a cell containing a preassigned spot originating from a finite universe. Thus, for any designated community induced by a finite universe, we only need to choose at least one spot belonging to some cell and appeal to the sizes of the underlying set to verify the union close set conjecture.
\end{remark}
\bigskip

\section{Main result}\label{sec:main}
In this section, we verify the union close set conjecture using the lemma \ref{frankil conjecture proof}.

\begin{theorem}[The Density law]\label{union close set conjecture proof}
Let $\mathbb{U}$ be a finite universe with an induced community $\mathcal{M}_{\mathbb{U}}$. There exists some spot $a_i\in \mathbb{U}$ such that 
\begin{align}
\mathcal{D}_{\mathcal{M}_{\mathbb{U}}}(a_i)\geq \frac{1}{2}.\nonumber
\end{align}
\end{theorem}

\begin{proof}
Let $\mathbb{U}$ be a finite universe with an induced community $\mathcal{M}_{\mathbb{U}}$. There exists some spot $a_i\in \mathbb{U}$ contained in some cell $\mathbb{A}\in \mathcal{M}_{\mathbb{U}}$. By appealing to Lemma \ref{frankil conjecture proof}, there exists some $l\geq 1$ such that 
$$
|\mathcal{M}_{\mathbb{U}}|\leq 2^{l}-1
$$ 
and 
\begin{align}
\# \left \{\mathbb{A}\in \mathcal{M}_{\mathbb{U}}|~a_i\in \mathbb{A}\right \}\geq 2^{l-1}\nonumber
\end{align}
so that we have the lower bound
\begin{align}
\frac{\# \left \{\mathbb{A}\in \mathcal{M}_{\mathbb{U}}|~a_i\in \mathbb{A}\right \}}{|\mathcal{M}_{\mathbb{U}}|}&\geq \frac{2^{l-1}}{2^l-1}\nonumber \\&=\frac{1}{2}\bigg(\frac{1}{1-\frac{1}{2^l}}\bigg).\nonumber
\end{align}
Taking the limits on both sides as $l\longrightarrow \infty$, the result follows immediately.
\end{proof}
\bigskip

It follows that for any finite universe $\mathbb{U}$ with an arbitrary induced community $\mathcal{M}_{\mathbb{U}}$, there exists some spot $a_i\in \mathbb{U}$ for which we can write the lower bound 
\begin{align}
\frac{\# \left \{\mathbb{A}\in \mathcal{M}_{\mathbb{U}}|~a_i\in \mathbb{A}\right \}}{|\mathcal{M}_{\mathbb{U}}|}&\geq \frac{1}{2}\bigg(\frac{1}{1-o(1)}\bigg).\nonumber
\end{align}

%%%%%%%%%%%%%%%%%%%%%%%%%%%%%%%%%%%%%%%%%%%%%%%%%%%%%%%%%%%%%%%%%%%%%%%%
\rule{100pt}{1pt}

\bibliographystyle{amsplain}

\end{document}